\documentclass[a4paper,12p]{article}
\title{Construction of Fractal Surfaces by Recurrent Fractal Interpolation Curves}
\author{\textrm{Chol-hui Yun, Hyong-chol O, Hui-chol Choi } $~~$ $~~$\\  \\
                  {\small\textit{Faculty of Mathematics}},
                        {\small\textit{\textbf{Kim Il Sung} University, Pyongyang, D. P. R. Korea}}\\   
}
\date{}  
\usepackage{latexsym}   
\usepackage[pdftex]{graphicx}
\usepackage{amsmath}
\usepackage{hyperref}

\begin{document}
\maketitle      
\begin{abstract}   
A method to construct fractal surfaces by recurrent fractal curves is provided. First we construct fractal interpolation curves using a recurrent iterated functions system (RIFS) with function scaling factors and estimate their box-counting dimension. Then we present a method of construction of wider class of fractal surfaces by fractal curves and Lipschitz functions and calculate the box-counting dimension of the constructed surfaces. Finally, we combine both methods to have more flexible constructions of fractal surfaces.
\end{abstract}
%
%
\section{Introduction}
A Fractal surface (or fractal curve) is a fractal set which is a graph of some continuous function on $\mathbf R^2$ (or $\mathbf R$). The method of  construction of fractal surfaces (or fractal curves) are closely related to the generation of fractal interpolation functions (FIF). The FIFs were introduced by Barnsley \cite{ban2} in 1986 and after that have been widely studied and used in approximation theory, image compression, computer graphics and modeling of natural surfaces such as rocks, metals, planets, terrains and so on. \\ \indent
The constructions of fractal surfaces, such as self-similar, self-affine or non self-affine surfaces, by IFSs or RIFSs have been studied in many papers (see \cite{BDD1, BD1, feng, mal, mas, MY, zha}). These surfaces are all attractors of some IFSs or RIFSs. In order to ensure the continuity of the surface, some authors in the early days assumed that the interpolation nodes on the boundary are collinear. Later Malysz \cite{mal}, Metzler and Yun \cite{MY} and Feng et al. \cite{feng} studied the construction of fractal interpolation surfaces on arbitrary data set using IFS. Furthermore Malysz \cite{mal} and Metzler and Yun \cite{MY} estimated the dimensions of the result surfaces. In \cite{BDD1}, the authors studied the construction of recurrent fractal interpolation surfaces.  In \cite{BDD1, mal} they use constant contraction factors in construction of IFS or RIFS and in \cite{feng, MY} they use function contraction factors in construction of IFS.\\\indent
The fractal properties of such rough surfaces as those of metals or rocks may expressed by sectional profiles of those surfaces. That's why constructions of fractal surfaces by fractal curves have been studied in many papers (see \cite{BD2, fal, man, WX, XFC,ZQL}). In \cite{man}, Mandelbrot suggested that the fractal dimension of a surfaces constructed by a single curve can be obtained by adding 1 to the fractal dimension of the curve. In \cite{fal}, Falconer introduced fractal surfaces constructed by the movement of a fractal curve along a segment and determined their fractal dimension. Xie and Fang proposed the so called star product fractal surfaces which are constructed by the movement of a fractal curves along another one (\cite{XFC}). In \cite{ZQL} a construction of fractal surfaces by 4 fractal curves, which are boundary curves of the constructed surface, and in \cite{WX} a construction of Bush type fractal surfaces by two Bush curves were studied. Bouboulis and Dalla \cite{BD2} provided a construction of fractal surfaces by the fractal interpolation. All the constructions of fractal surfaces by fractal curves have common property that the fractal dimension of the constructed surfaces is determined by one of the fractal curves constructing them, and that the fractal curves are combined with constants or one dimensional functions, which are called combination functions. And in \cite{ban1, ban3}, they constructed the fractal curves by IFS or RIFS with constant vertical scaling factors. These constructions constrain vertical scaling factors (to be constants) and combination functions (to be constants or one-dimensional functions).
 Hence, these constructions lack the flexibility needed to model complex natural curves and surfaces.
Because objects in nature modeled by fractals are very complex and irregular, we need flexible constructions which can construct broader class of fractal sets for modeling them. The aim of this paper is to introduce constructions of a new, broader class of fractal curves and surfaces.\\\indent
The recurrent iterated function system (RIFS) is a generalization of the IFS. The vertical scaling factors of transformations in IFS (or RIFS) are contraction ratios at each point of an attractor of IFS (or RIFS) and characterize a fractal structure of the attractor. And it is more general that each point of images in nature has a different contraction ratio, and the class of Lipschitz functions is broader than the one of constants and the one-dimensional functions. A fractal dimension is important parameter that characterizes a fractal structure of set. \\\indent
 In order to ensure more flexibility in modeling natural shapes and phenomena or in image processing we introduce a construction of recurrent fractal interpolation curves by more general RIFSs with function scaling factors (Theorem 1, which ensures that an attractor of constructed RIFS is a graph of some continuous function and gives how to make the fractal interpolation function, which can be used to get interpolation values from a discrete data set with fractal property) and estimate the box-counting dimension (Theorem 2, which gives the fractal structure of the constructed curves). And we construct broader class of fractal surfaces by fractal curves and Lipschitz functions, and calculate the box-counting dimension of the constructed surfaces (Theorem 3, which shows the b-counting dimension of surfaces constructed by the recurrent fractal and Lipschitz functions and shows that the constructed surfaces are fractal surfaces, and their fractal structure, and that the box-counting dimension of the constructed fractal surfaces can be controlled by one of the recurrent fractal curves). Finally we combine both constructions to construct fractal surfaces. These constructions of recurrent fractal interpolation curves and fractal surfaces can be used to approximate discrete sequence of data (like seismic, electrocardiogram data), natural curves (like coastlines, profiles of mountain ranges, tops of clouds and horizons over forest) and natural surfaces (like rocks, metals, terrains), respectively.    \\\indent
The remainder of the article is organized as follows: The section 2 describes some basic notions of RIFS. The section 3 constructs recurrent fractal curves with function vertical scaling factors and estimates their box counting dimensions. The section 4.1 provides a general method of construction of fractal surfaces combining fractal curves using Lipschitz functions and a formula of the box counting dimensions. As one application of the results of the sections 3 and 4.1, the section 4.2 describes a method using some Lipschitz functions and provides some examples.  
%
\section{Preliminaries on RIFS}
In this section we describe some basic notions on recurrent iterated function systems and related lemmas. The class of irreducible matrices are important in RIFS theory.

\textbf{Definition 1}\cite{Gan} Let $Q=(q_{ij})$ be $n\times n$matrix. $Q$ is called \textit{\textbf{reducible}}, if the set of indices $\{1,2, \ldots, n\}$ can be decomposed into two distinct and complementary subsets $\{i_1, i_2,\ldots,i_{\mu } \}$ and $\{j_1, j_2,\ldots,j_{\nu }\}$, i.e.,
\begin{equation}
\{ i_1, i_2,\ldots , i_{\mu }\} \cap \{ j_1, j_2,\ldots , j_{\nu }\} = \emptyset , \nonumber
\end{equation}
\begin{equation}
\{ i_1, i_2,\ldots , i_{\mu }\} \cup  \{ j_1, j_2,\ldots , j_{\nu } \} = \{ 1, 2, \ldots , n \} \nonumber
\end{equation}
such that 
\begin{equation}
q_{i_a,j_b} = 0, a=1,2,\ldots ,\mu ; b=1,2, \ldots , \nu . \nonumber 
\end{equation}
If this is not possible, then $Q$ is an \textit{\textbf{irreducible}} matrix.

The following results on irreducible matrices are useful and can refer to \cite{Gan,se}.
\newtheorem{lem}{Lemma}
\begin{lem}\label{L:001}
Consider an square matrix $A$. The following statements are equivalent. \\\indent
\textnormal{(1)} $A$ is irreducible. \\\indent
\textnormal{(2)} For every pair $(i,j), i,j=1,2,\ldots,n$, there is a $k>0$ such that the element of the $i-$th row $j-$th column of the matrix $A^k$ is positive.
\end{lem}
\begin{lem}\label{L:002}
$A$ is a non-negative, reducible $n\times n$ matrix, iff the matrix $(I_n+A)^{n-1}$ is positive (where $I_n$ is the identity matrix). 
\end{lem}

\begin{lem}\label{L:003}
(\textbf{Perron-Frobenius Theorem}) Let $A \geq 0$ be an irreducible square matrix. Then we have the following two statements.\\\indent
\textnormal{(1)} The spectral radius $\rho (A)$ of $A$ is an eigenvalue of $A$ and it has strictly positive eigenvector $y$ (i.e., $y_i>0$ for all $i$)\\\indent
\textnormal{(2)} $\rho (A)$ increases if any element of $A$ increases.
\end{lem}

\textbf{Definition 2} \cite{ban3, Bou} $(X, d_X)$ is a complete metric space and $X_1, X_2, \ldots, X_n$ are subsets of $X$. Let $W_i:X_i \rightarrow X(i=1,\ldots,n)$ be contraction maps. $M=(p_{ij})^n_{i,j=1}$ is a $n \times n$ row-stochastic matrix (i.e., such a matrix that $\sum^n_{j=1}p_{ij}=1, i=1,\ldots ,n$) and irreducible. Then $\left\{X;M;W_1,\ldots,W_n\right\}$ is called a \textit{\textbf{recurrent iterated function system}}(or simply RIFS). The $n\times n$ matrix $C=(c_{ij})$ defined by
\begin{equation}\label{c02}
c_{ij}=\left\{
\begin{array}{ll}
1, & p_{ji}>0,\nonumber\\
0, & p_{ji}=0.
\end{array}
\right.
\end{equation}
is called the \textit{\textbf{connection matrix}} of RIFS.\\\indent

Let $\mathcal{H}(X)$ be  the set of all nonempty compact subset of $X$ and $h$ the Hausdorff distance in $\mathcal{H}(X)$. Then  $(\mathcal{H}(X),h)$ is a complete metric space \cite{ban1}. In the set $\tilde{\mathcal{H}}(X):=\underbrace{\mathcal{H}(X)^n=\mathcal{H}(X) \times \dots \times \mathcal{H}(X)}_n$, we define a distance $\tilde{h}:\tilde{\mathcal{H}}(X) \times \tilde{\mathcal{H}}(X)\rightarrow  \mathbf{R}$ as follows:
\begin{eqnarray}
\tilde{h}((A_1,A_2,\dots ,A_n),(B_1,B_2,\dots ,B_n)):=\max_{i=1,\dots ,n} h(A_i,B_i), \nonumber \\
\forall (A_1,A_2,\dots ,A_n),(B_1,B_2,\dots ,B_n)\in \tilde{\mathcal{H}}(X) \nonumber 
\end{eqnarray}
Then $(\tilde{\mathcal{H}}(X), \tilde{h})$ is a complete metric space \cite{Bou}.

For a given recurrent iterated functions system $\left\{X;M;W_1,\ldots,W_n\right\}$, we define the transformation $\mathbf{W}:\tilde{\mathcal{H}}(X)\rightarrow \tilde{\mathcal{H}}(X)$ of $\tilde{\mathcal{H}}(X)$ as follws: let

$$
\mathbf{W}(\mathbf{B}):=\left(
\begin{array}{ccc}
c_{11}W_1(B_1)\cup c_{12}W_1(B_2)\cup \dots \cup c_{1n}W_1(B_n) \\
c_{21}W_2(B_1)\cup c_{22}W_2(B_2)\cup\dots\cup c_{2n}W_2(B_n)\\
\vdots\\
c_{n1}W_n(B_1)\cup c_{n2}W_n(B_2)\cup\dots\cup c_{nn}W_n(B_n)
\end{array}
\right)
=\left(
\begin{array}{ccc}
\bigcup_{j\in\Lambda(1)}W_1(B_j)  \\
\bigcup_{j\in\Lambda(2)}W_2(B_j)  \\
\vdots\\
\bigcup_{j\in\Lambda(n)}W_n(B_j)  \\
\end{array}
\right)
$$
for every $\mathbf{B}=(B_1,\dots, B_n)\in \tilde{\mathcal{H}}(X)$. Here $\Lambda(i)=\left\{j:c_{ij}=1\right\},i=1,\dots,n$. We will denote the transformation $\mathbf{W} $ as a matrix form by
$$
\mathbf{W}=\left(
\begin{array}{cccc}
c_{11}W_1 & c_{12}W_1 & \dots & c_{1n}W_1\\
c_{21}W_2 & c_{22}W_2 & \dots & c_{2n}W_2\\
   \vdots & \vdots    & \ddots & \vdots\\
c_{n1}W_n & c_{n2}W_n & \dots & c_{nn}W_n
\end{array}
\right)
$$

This transformation $\mathbf{W} $ is a contraction map in $\tilde{\mathcal{H}}(X)$ and the unique fixed point $\mathbf{A}=(A_1,\dots, A_n) $ is called the \textit{\textbf{attractor}} or \textit{\textbf{invariant set}} of the RIFS, or \textit{\textbf{recurrent fractal}} \cite{ban3, Bou}. \\\indent
\textbf{Note}. The invariant set of RIFS is a vector whose elements are nonempty compact sets. If $\mathbf{A}=(A_1,\dots, A_n) $  is the invariant set of  a RIFS$\left\{X;M;W_1,\ldots,W_n\right\}$, then we have $A_i=\bigcup_{j\in\Lambda(i)}W_i(A_j)$, $ \forall i=1,\dots,n$. Usually, making a slight abuse of notation, we often call the union $ A=\cup_{i=1}^nA_i$ of all $A_i$ as the attractor of RIFS, too.

\textbf{Definition 3}. A curve which is the invariant set of a RIFS on $\mathbf{R}^2$ is called a \textit{\textbf{recurrent fractal curve}}(or RFC).

\textbf{Definition 4}. In general, the \textit{\textbf{box-counting dimension}} $\mathrm{dim}_\mathrm{B} \mathcal{A}$ of a fractal set $\mathcal{A}$ is defined by
\begin{displaymath}
\mathrm{dim_{B}}\mathcal{A}=\lim_{\delta \to 0}\frac{\log N_{\delta}\left(\mathcal{A}\right)}{-\log \delta}
\end{displaymath}
(if this limit exists), where $N_{\delta}(\mathcal{A})$ is any of the followings (see \cite{feng}):
\begin{enumerate}
\item[(i)] the smallest number of closed balls of radius $\delta$ that cover $\mathcal{A}$; 
\item[(ii)] the smallest number of cubes of side $\delta$ that cover $\mathcal{A}$; 
\item[(iii)] the number of $\delta$-mesh cubes that intersect $\mathcal{A}$; 
\item[(iv)] the smallest number of sets of diameter at most $\delta$ that cover $\mathcal{A}$; 
\item[(v)] the largest number of disjoint balls of radius $\delta$ with centers in $\mathcal{A}$.
\end{enumerate}
Here we use (iii).
%
%
\section{Construction of recurrent fractal interpolation curves}
In this section, we construct recurrent fractal curves with function vertical scaling factors and estimate their box-counting dimension.\\\indent
Let a data set be 
$
P=\left\{(x_i,y_i)\in \mathbf{R}^2;i=1,\ldots,n\right\}, (x_0<x_1<\ldots <x_n) 
$
and let
$
N_n=\left\{1,\ldots,n\right\}, I=[x_0, x_n], I_i=[x_{i-1}.x_{i}], i=1,\ldots ,n .
$
We denote  Lipschitz (or contraction) constant of Lipschitz (or contraction) mapping $f$  by $L_f$ (or $c_f$). Let $l\geq 2, l\in \mathbf{N}$ and let $\tilde{I}_k=[x_{s(k)},~x_{e(k)}],~x_{s(k)}$,$~x_{e(k)}~\in \left\{x_0,\ldots,x_n\right\}$, here $e(k)-s(k) \geq 2,~k=1,\ldots,l $, as a custom, $I_i$ is called a region and $\tilde{I}_k$ a domain.  \\\indent
We collect a map $\gamma :N_n\rightarrow N_l$. This means that we relate every region to a domain. For each $i\in N_n$, let $k=\gamma(i)$. For $i\in N_n$(and $k=\gamma (i)$), let a mapping $L_i=L_{i,k}:\tilde{I}_k\rightarrow I_i$ be such a contraction homeomorphism that $L_{i,k}:\left\{x_{s(k)},x_{e(k)}\right\}\rightarrow \left\{x_{i-1},x_i\right\}$. Such a map can be easily constructed in a standard way. Let $H \subset \mathbf{R}$ be such a sufficiently large interval that $y_i \in H, i=1,\dots,n$. Let $F_i=F_{i,k}:\tilde {I}_k\times H\rightarrow \mathbf{R}$ be a function defined by  $F_{i,k}(x,y)=s_{i,k}(L_{i,k}(x))a(y)+b_{i,k}(x)$  and satisfy 
\begin{eqnarray}\label{e1} 
F_{i,k}(x_{\alpha},y_{\alpha})=y_{\beta}, \alpha\in \left\{s(k),e(k)\right\},
\end{eqnarray}
where  $L_{i,k}(x_{\alpha})=x_{\beta}, \beta \in \left\{i-1,i\right\}$ and $a(y)$ is a Lipschitz function on $H$ such that $a(y_\alpha)=y_\alpha, \alpha\in \left\{s(k),e(k)\right\}$, $b_{i,k}(x)$  is a Lipschitz function defined on the domain $\tilde{I}_k $ and $s_{i,k}(x)$ is a contraction function on $I_i$  with $|s_{i,k}(x)L_a|<1$ (which is called a \textit{\textbf{vertical scaling factor}}.)\\\indent

\textbf{Example 1} The function $F_{i,k}$ satisfying (\ref{e1}) can be easily constructed. When $h(x), g(x)$  are Lipschitz mappings on $I$ and satisfy the  conditions
 $$
 g(x_\alpha)=y_{\alpha}, \alpha\in\left\{s(k),e(k)\right\},
 h(x_i)=y_i,i=0,1,\ldots,n,
 $$
and $s_{i,k}(x)$  are taken as free unknown functions, the functions
$$
F_{i,k}(x,y)=s_{i,k}(L_{i,k}(x))(a(y)-g(x))+h(L_{i,k}(x)),
$$
satisfy (\ref{e1}). 

We define transformations  $W_{i,k}:\tilde{I}_k\times H \rightarrow I_i\times \mathbf{R}~(i=1,\ldots,n;~k=\gamma(i))$ by
\begin{equation*}
W_{i, k}(x,~y)=(L_{i,k}(x),~F_{i,k}(x,~y)).
\end{equation*}
Then there exists some distance equivalent to the Euclidean metric on  $\mathbf{R}^2$ such that  $W_{i,k}(i=1,\ldots,n; ~k=\gamma(i))$ are contraction transformations with respect to the distance.\\\indent
We define a \textit{\textbf{row-stochastic matrix}}  $M=(p_{ij})_{n\times n}$  by 
 \begin{equation}\label{c1}
p_{ij}=\left\{
\begin{array}{ll}
\frac{1}{a_i}, & I_i\subseteq \tilde{I}_{\gamma (j)}\nonumber\\
0 ,& \mathrm{otherwise}
\end{array}\right.,
\end{equation}
where for every $i(\in N_n)$,  the number $a_i$  indicates the number of the domains $\tilde{I}_{\gamma (j)}$ containing the region $I_i$. Then a \textit {\textbf{connection matrix}} $C=(c_{ij})_{n \times n}$  is defined as follows
 \begin{equation}\label{c2}
c_{ij}=\left\{
\begin{array}{ll}
1, & p_{ji}>0,\nonumber\\
0, & p_{ji}=0.
\end{array}\right.
\end{equation}
It is clear that if  the row-stochastic matrix is irreducible, then the connection matrix is also irreducible.\\\indent
Let $E=I\times H \subset  \mathbf{R}^2$. We denote the attractor of the recurrent iterated functions system (RIFS) $\{ E; M;W_{i,k},i=1,\ldots,n,~\gamma(i)=k\in\left\{1,\ldots,l\right\}\}$ by $\mathcal{A}$. Then the following theorem shows that $\mathcal{A}$   is a recurrent fractal curve.\\\indent
\newtheorem{Thot}{Theorem}
\begin{Thot}\label{T:1}
The attractor $\mathcal{A}$   constructed above is a graph of some continuous function which interpolates the data set $P$.
\end{Thot}
\textbf{Proof.} Let  $\mathrm{C}(I)=\left\{\varphi \in C^0(I);\varphi (x_i)=y_i, i=0,~1,~\ldots~,~n \right\}$, then the set $\mathrm{C}(I)$  is complete metric space with respect to norm $||\cdot ||_\infty $. We can easily know that the operator $T:\mathrm{C}(I)\rightarrow \mathrm{C}(I)$ ; $ (T\varphi )(x)=F_{i,k}(L^{-1}_{i,k}(x),\varphi (L^{-1}_{i,k}(x))), x\in I_i$ is well defined and the operator  $T$ is a contraction on the complete metric space $\mathrm{C}(I)$. Therefore the operator  $T$ has a unique fixed point in $\mathrm{C}(I)$, which we denote by $f$. Then the $f$ is presented by
 $$f(x)=s_{i,k}(x)f(L^{-1}_{i,k}(x)+b(x), i=1,\ldots,n,  \gamma(i)=k$$
which means that $Gr(f)=\mathcal{A}$ , where $Gr(f)$  denote a graph of $f$.$~~~~~\Box $\\\indent
We estimate box-counting dimensions of the recurrent fractal curves constructed above. We can assume that $I=[0,1]$, since the box-counting dimension is invariant under bi-Lipschitz mapping.\\
For a set $D(\subset \mathbf{R}^1$ or $\mathbf{R}^2)$ and a function $f$ defined on $D$, 
\begin{displaymath}
R_f[D]=sup\left\{|f(x_2)-f(x_1)|;x_1,x_2\in D\right\}.
\end{displaymath}
is called the \textit{maximum variation} of $f$ on $D$.\\\indent
Let  $I$ be an interval in $\mathbf{R}$, $L:I\rightarrow I$  a contraction homeomorphism, $a,b:\mathbf{R}\rightarrow \mathbf{R}$  Lipschitz mappings and  $s:I\rightarrow \mathbf{R}$  contraction mappings with $|s(x)La|<1$ . We define a mapping  $F:I\times \mathbf{R}\rightarrow \mathbf{R}$ by
 \begin{displaymath}
F(x,y)=s(L(x))a(y)+b(x), (x,y)\in I\times \mathbf{R}.
\end{displaymath}
\begin{lem}\label{L:1}
 Let $f:I\rightarrow  \mathbf{R}$  be continuous function. Then
 \begin{displaymath}
 R_{F(L^{-1},~f\circ L^{-1})}[L(I)]\leq\bar s L_a R_f[I]+|I|(c_s\bar a_f+L_b).
 \end{displaymath}
Here $|I|$ is a length of the interval $I$,  $\bar s=\max_I\left|s(x)\right|~$ and $~\bar a_f=\max_I\left|a(f(x))\right|$  .
\end{lem}
\textbf{Proof.} For any $x(\in L(I))$, let denote $\tilde x=L^{-1}(x)(\in I)$. Then 
\begin{eqnarray}
&&|F(L^{-1}, f\circ L^{-1})(x)-F(L^{-1}, f\circ L^{-1})(x')|=\nonumber \\
&&~~~~=|F(L^{-1}(x), f\circ L^{-1}(x))-F(L^{-1}(x'), f\circ L^{-1}(x'))|\nonumber \\
&&~~~~=|s(x)a(f(\tilde x))+b(\tilde x)-s(x')a(f(\tilde x'))-b(\tilde x')|\nonumber \\
&&~~~~=|s(x)a(f(\tilde x))-s(x)a(f(\tilde x'))+s(x)a(f(\tilde x'))-s(x')a(f(\tilde x'))+b(\tilde x)-b(\tilde x')|\nonumber \\
&&~~~~\leq\bar s L_a R_f[I]+c_s|x-x'|\bar a_f+L_b |\tilde x-\tilde x'|\nonumber \\
&&~~~~\leq\bar s L_a R_f[I]+|I|(c_s\bar a_f+L_b).\nonumber 
\end{eqnarray}
$\Box $\\   \indent
We assume that the row-stochastic matrix $M$ is irreducible and the mapping $a(x)$ is identity. Note that in what follows the latter \textquoteleft $a$' is used in another meaning. 
Let $x_{i+1}-x_i=\frac{1}{n}(i=0,1,\ldots,n-1), x_{e(k)}-x_{s(k)}=\frac{a}{n}(k=1,\ldots,l; a\in \mathrm{N})$. Let $L_{i,k}(i\in N_n, k=\gamma(i))$ be similitude contraction transformations. Then the number of $I_j$ contained in $\tilde{I}_k$  is $a$ . Let $\bar S$ and $\underline{S}$ be diagonal matrices  
$$ \bar S=diag(\bar s_1,\ldots, \bar s_n),~~\underline S=diag(\underline s_1,\ldots, \underline s_n),$$
respectively, where $\bar s_i=max_{I_i}|s_{i,\gamma(i)}(x)|$, $\underline s_i=min_{I_i}|s_{i,\gamma(i)}(x)|$. 
%
\begin{Thot}\label{T:2}
Let $\mathcal{A}$ be the recurrent fractal curve in Theorem~\ref{T:1}. If there is some interval $\tilde{I}_{k_0}$  such that the points of $\mathrm{P}\cap (\tilde{I}_{k_0}\times\mathbf{R})$  are not collinear, then the box-counting dimension $\mathrm{dim}_B\mathcal{A}$  of $\mathcal{A}$  has the following lower and upper bounds;\\\indent
\textnormal{(1)} If $\underline{\lambda}\leq 1$, then $$1+\log_a\underline{\lambda}~\leq \mathrm{dim}_B\mathcal{A} \leq 1+\log _a\bar{\lambda}~,$$\indent
\textnormal{(2)} If $\bar{\lambda}\leq 1$, then $$\mathrm{dim}_B\mathcal{A}=1~,$$
where $\underline{\lambda}=\rho(\underline{S} \mathrm{C})$ and $\bar\lambda=\rho(\bar S \mathrm{C})$ are respectively spectral radii of the irreducible matrices $\underline{S} \mathrm{C}$ and $\bar S \mathrm{C}$.
\end{Thot}
\textbf{Proof.}  \textbf{Proof of (1)}: Let $f$ be a contraction function whose graph is $\mathcal{A}$. We denote $R_f[\tilde{I}_k]$  by $R_k$  and  $\frac{1}{a^r}$  by $\varepsilon _r$ for simplicity. Then $r\rightarrow\infty  \Longleftrightarrow \varepsilon _r\rightarrow 0$ .\\
After applying each $W_{i,k}$ to the interpolation points in the interval $\tilde{I}_k$ one time, we have $(a+1)$ new points in every interval $I_i(i \in N_n)$. According to the hypothesis, the interpolation points lying inside $\tilde{I}_{k_0}$ are not collinear and the connection matrix $C$  is irreducible. Thus for every interval $I_i$ there are three points which are not collinear, the maximum vertical distance (computed only along the $y$-axis) from one of the three points to the line through other two interpolation points is greater than 0. The maximum value is called a \textit{height} and denoted by $H_i$ . \\\indent
By Lemma~\ref{L:1}, on each region $I_i$ we have
$$ R_f[I_i]\leq \bar {s}_i~R_k+\frac{a}{n}e, $$
where $e=c_s\bar f+L_b$.\\\indent
We define non negative vectors $\mathbf{h}_1, \mathbf{r}, \mathbf{u}_1$ and $\mathbf{i}$  by
 \begin{equation}
\mathbf{h_1}=\left[
     \begin{array}{ll}
               H_1\nonumber\\
                ~\cdot \nonumber\\
                ~\cdot \nonumber\\
                ~\cdot \nonumber\\
                H_n \nonumber
   \end{array}
   \right ], ~~~\mathbf{r}=\left[
                \begin{array}{ll}
                     \bar s_1R_1\nonumber\\
                     ~~\cdot \nonumber\\
                     ~~\cdot \nonumber\\
                     ~~\cdot \nonumber\\
                     \bar s_nR_n \nonumber
                \end{array}
           \right ],~~~\mathbf{i}=\left[
\begin{array}{ll}
1\nonumber\\
\cdot \nonumber\\
\cdot \nonumber\\
\cdot \nonumber\\
1 \nonumber
\end{array}
\right], ~~~\mathbf{u_1}=\mathbf{r}+\frac{a}{n}e\mathbf{i}.
\end{equation}
Since  $\mathcal{A}$ is the graph of a continuous function defined on $I$, the smallest number of  $\varepsilon _r$-mesh squares  necessary to cover $I_i\times\mathbf{R}\cap \mathcal{A}$  is greater than the smallest number of  $\varepsilon _r$-mesh squares   necessary to cover vertical line with the length $H_i$ and less than the smallest number of  $\varepsilon _r$-mesh squares necessary to cover a rectangle $I_i\times[\underline{f}_i,~\bar f_i]$, where $$\underline{f}_i=\min_{I_i}|f(x,~y)|,~~\bar f_i=\max_{I_i}|f(x,~y)|.$$
	Therefore,  
\begin{displaymath} 
\sum^n_{i=1}[H_i\varepsilon ^{-1}_r]-n\leq N(\varepsilon_r)\leq \sum^n_{i=1}\left(\left[\left(\bar s_iR_i+\frac{a}{n}e\mathbf{i}\right)\varepsilon^{-1}_r\right]+1\right)\left(\left[\frac{\varepsilon^{-1}_r}{n}\right]+1\right),
\end{displaymath}
i.e.
\begin{displaymath} 
\Phi(\mathbf{h_1}\varepsilon ^{-1}_r)-n\leq N(\varepsilon_r)\leq \Phi(\mathbf{u_1}\varepsilon ^{-1}_r+\mathbf{i})\left(\left[\frac{\varepsilon^{-1}_r}{n}\right]+1\right),
\end{displaymath}
where $\frac{1}{\mathrm{n}}>\varepsilon_r$ and 
$$\Phi(\mathbf{a})=a_1+\dots+a_n, ~~\mathbf{a}=(a_1,\ldots,a_n).$$\indent
After applying  $W_{i,k}$ twice, we have $a$  subintervals of length $\frac{1}{an}$  in each $I_i(i\in N_n)$. Those subintervals are mapped by the transformation $W_{i,k}$  from subintervals lying inside the intervals $\tilde{I}_{\gamma(i)}$ corresponding to the interval $I_i$ containing themselves, and thus for each subinterval $I_i$  the height on the new subintervals produced in $I_i$ is not less than $\underline{s}_i\cdot H$, where $H$ is the height on original subinterval contained in the interval $\tilde{I}_{\gamma(i)}$. Therefore, the sum of maximum variance of $f$  on $a$  subintervals of the length $\frac{1}{an}$  contained in the interval $I_i(i\in N_n)$  is not greater than $i$-th coordinate of vector
$ \mathbf{u}_2=\bar{S}C\mathbf{u}_1+\frac{a}{n}e\mathbf{i}$, the sum of the heights is not less than  $i$-th coordinate of vector $ \mathbf{h}_2=\underline{S}C\mathbf{h}_1$ and
\begin{displaymath}
\mathbf{\Phi}(\mathbf{h}_2\varepsilon^{-1}_r)-an\leq N(\varepsilon_r)\leq \mathbf{\Phi}(\mathbf{u}_2\varepsilon^{-1}_r+a\mathbf{i})\left(\left[\frac{\varepsilon^{-1}_r}{an}\right]+1\right),
\end{displaymath}
where $\frac{1}{an}>\varepsilon_r$. \\\indent
By induction  after taking $k$  such that $a\varepsilon_r \geq \frac{1}{a^{k-1}n}\geq \varepsilon_r$, that is, 
$r-\log_an+1>k\geq r-\log_an$ and applying $W_{j,\gamma(i)}$  $k$ times, we get  $a^{k-1}$  subintervals of the length $\frac{1}{a^{k-1}n}$  contained in each interval $I_i$ and
\begin{eqnarray}\label{eq:2}
\mathbf{\Phi} (\mathbf{h}_k\varepsilon^{-1} _r)-a^{k-1}n \leq N(\varepsilon _r)\leq \mathbf{\Phi}(\mathbf{u}_k\varepsilon^{-1} _r+a^{k-1}\mathbf{i})\left(\left[\frac{\varepsilon^{-1} _r}{a^{k-1}n}\right]+1\right),
\end{eqnarray}
where $\mathbf{u}_k=\bar S\mathrm{C}\mathbf{u}_{k-1}+\frac{a^k}{n}b\mathbf{i},~\mathbf{h}_k=\underline{S}\mathrm{C}\mathbf{h}_{k-1}$. Then we have
\begin{eqnarray*}		
&&\mathbf{u}_k=(\bar SC)^{k-1}\mathbf{r}+(\bar SC)^{k-1}\frac{a}{n}e\mathbf{i}+(\bar SC)^{k-2}\frac{a}{n}e\mathbf{i}+\ldots+(\bar SC)\frac{a}{n}e\mathbf{i}+\frac{a}{n}e\mathbf{i},\\
&&\mathbf{h}_k=(\underline{S}C)^{(k-1)}\mathbf{h}_1.	
\end{eqnarray*}
Since $\underline{\textit S}C$ and $\bar{S}C$ are non-negative irreducible matrix, from Perron-Frobenius's theorem (Lemma \ref{L:003} ) there exist strictly positive eigenvectors of $\underline{\textit S}C$ and $\bar{S}C$(which correspond to eigenvalues $\underline{\lambda}=\rho(\underline{\textit S}C)$ and $\bar{\lambda}=\rho(\bar{S}C)$ of $\underline{\textit S}C$ and $\bar{S}C$) and we can choose such strictly positive eigenvectors $\bar{\textbf e}$ , $\underline{\textbf e}$ (which correspond to eigenvalues $\underline{\lambda}$ , $\bar{\lambda}$ respectively) that 
 $$\mathbf{r}\leq \bar{\textbf e}, ~~~~b\mathbf{i}<n\bar{\textbf e}, ~~~~0<\underline{\textbf e}_1<\mathbf{h}_1~.$$
Then 
\begin{eqnarray}
&& N(\varepsilon _r)\leq \mathbf{\Phi}(\mathbf{u}_k\varepsilon^{-1} _r+a^{k-1}\mathbf{i})\left(\left[\frac{\varepsilon^{-1} _r}{a^{k-1}n}\right]+1\right)\nonumber\\
&&~~~\leq \mathbf{\Phi}(\mathbf{u}_k\varepsilon^{-1} _r+a^{k-1}\mathbf{i})(a+1)\nonumber \\
&&~~~\leq \mathbf{\Phi}( (\bar SC)^{k-1}\mathbf{r}\varepsilon^{-1}_r+(\bar SC)^{k-1}\frac{a}{n}e\mathbf{i}\varepsilon^{-1} _r+(\bar SC)^{k-2}\frac{a}{n}e\mathbf{i}\varepsilon^{-1} _r+\ldots+(\bar SC)\frac{a}{n}e\mathbf{i}\varepsilon^{-1} _r+\nonumber\\
&&~~~~~~~+\frac{a}{n}e\mathbf{i}\varepsilon^{-1} _r+a^{k-1}\mathbf{i})(a+1) \nonumber \\
&&~~~\leq \mathbf{\Phi}( (\bar SC)^{k-1}\mathbf{\bar e}\varepsilon^{-1} _r+(\bar SC)^{k-1}\mathbf{\bar e}\varepsilon^{-1} _r+(\bar SC)^{k-2}\bar{\textbf e}\varepsilon^{-1}_r+\ldots+(\bar SC)\bar{\textbf e}\varepsilon^{-1} _r+ \nonumber \\
&&~~~~~~~+\mathbf{\bar e}\varepsilon^{-1} _r+a^{k-1}\mathbf{i})(a+1) \nonumber \\
&&~~~=(\bar{\lambda}^{k-1}\Phi(\bar{\textbf e})\varepsilon^{-1}_r+\bar{\lambda}^{k-1}\Phi(\bar{\textbf e})\varepsilon^{-1}_r+\bar{\lambda}^{k-2}\Phi(\bar{\textbf
 e})\varepsilon^{-1}_r+\ldots+\bar{\lambda}\Phi(\bar{\textbf e})\varepsilon^{-1} _r+ \nonumber \\
&&~~~~~~~+\Phi(\bar{\textbf e})\varepsilon^{-1}_r+a^{k-1}\Phi(\mathbf{i})(a+1) \nonumber \\
&&~~~\leq (\bar{\lambda}^{r-\nu }\Phi(\bar{\textbf e})\varepsilon^{-1}_r+\bar{\lambda}^{r-\nu }\Phi(\bar{\textbf e})\varepsilon^{-1}_r+\bar{\lambda}^{r-\nu -1}\Phi(\bar{\textbf
 e})\varepsilon^{-1}_r+\ldots+\bar{\lambda}\Phi(\bar{\textbf e})\varepsilon^{-1} _r+ \nonumber \\
&&~~~~~~~+\Phi(\bar{\textbf e})\varepsilon^{-1}_r+a^{r-\nu }\Phi(\mathbf{i})(a+1), \nonumber
\end{eqnarray}
where  $\nu:=\mathbf{\log}_an$.\\\indent
On the other hand, since $(\underline{s}_i\leq \bar{s}_i$ for $i=1,\ldots,n$, from Perron-Frobenius's theorem (Lemma \ref{L:003}) we have $\underline{\lambda}\leq \bar{\lambda}$. Therefore if $\underline{\lambda}>1$ , then $1>\frac{1}{\underline{\lambda}}\geq \frac{1}{\bar{\lambda}}$ , we obtain
 \begin{eqnarray} 
&& N(\varepsilon _r)\leq \bar{\lambda}^{r-\nu}\Phi(\bar{\textbf e})\varepsilon^{-1}_r \left (1+1+\frac{1}{\bar{\lambda}}+\ldots+\frac{1}{\bar{\lambda}^{r-\nu}}+\frac{\mathrm{n}}{\bar {\lambda}^{r-\nu}\Phi(\bar{\textbf e})a^\nu}\right )(a+1)\nonumber \\
&&~~~~~~=\bar{\lambda}^r\varepsilon^{-1}_r\bar{\lambda}^{-\nu}\Phi(\bar{\textbf e}) \left (1+\frac{1-\left(\frac{1}{\bar{\lambda}}\right)^{r-\nu +1}}{1-\frac{1}{\bar{\lambda}}}+\frac{\mathrm{n}}{\bar{\lambda}^{r-\nu}\Phi(\bar{\textbf e})a^\nu }\right)(a+1).\nonumber 
\end{eqnarray}
Here let $\delta (r):=\bar{\lambda}^{-\nu}\Phi(\bar{\textbf e}) \left (1+\frac{1-\left(\frac{1}{\bar{\lambda}}\right)^{r-\nu +1}}{1-\frac{1}{\bar{\lambda}}}+\frac{\mathrm{n}}{\bar{\lambda}^{r-\nu}\Phi(\bar{\textbf e})a^\nu }\right)(a+1),$ then $\delta (r)>0$  and
\begin{eqnarray}\label{eq:3}
\dim_B\mathcal{A}=\lim_{\varepsilon_r\rightarrow 0}\frac{\log N(\varepsilon _r)}{-\log \varepsilon _r}\leq 1+\log_a\bar {\lambda}.	
\end{eqnarray}	
By (\ref{eq:2}), we have
\begin{eqnarray*}
&& N(\varepsilon_r)\geq \mathbf{\Phi}(\mathbf{h}_k\varepsilon^{-1}_r)-a^{k-1}{n}~~~~~~~=\mathbf{\Phi}((\underline{\textit S}C)^{k-1}\mathbf{h}_1\varepsilon^{-1} _r)-a^{k-1}{n}\\
&&~~~~~~ \geq \mathbf{\Phi}((\underline{\textit S}C)^{k-1}\mathbf{\underline{e}}\varepsilon^{-1}_r)-a^{k-1}{n}=\underline{\lambda}^{k-1}\Phi(\underline{\textbf e})\varepsilon^{-1}_r-a^{k-1}{n}\\
 &&~~~~~~ \geq \underline{\lambda}^{r-\nu-1}\Phi(\underline{\textbf e})\varepsilon^{-1}_r-a^{r-2\nu}n\varepsilon^{-1}_r\\
&&~~~~~~=\varepsilon^{-1}_r\underline{\lambda}^{r}\left(\underline{\lambda}^{-\nu-1}\Phi(\underline{\textbf e})-\frac{a^{-\nu}n}{\underline{\lambda}^{r}}\right). 
\end{eqnarray*}
Since $\underline{\lambda}>1$ , there is $r_0$ such that $\eta (r):=\underline{\lambda}^{-\nu-1}\Phi(\underline{\textbf e})-\frac{a^{-\nu}n}{\underline{\lambda}^{r}}>0$ for any $r(>r_0)$ and therefore we have 
\begin{equation}\label{eq:4}
\frac{\log N(\varepsilon _r)}{-\log \varepsilon _r}\geq 1+\mathbf{\log}_a{\underline{\lambda}}+\frac{1}{r}\log_a{\eta (r)},~~ r(>r_0)
\end{equation}  
By (\ref{eq:3}), (\ref{eq:4}) if $\underline{\lambda}>1$ , then we have
 \begin{displaymath}
1+\log_a{\underline{\lambda}}\leq  \dim_B\mathcal{A}\leq 1+\log_a{\bar{\lambda}}.
\end{displaymath}
\textbf{Proof of (2)}: If $\bar{\lambda}\leq 1$, we have  
\begin{eqnarray}	
&& N(\varepsilon _r)\leq (\bar{\lambda}^{r-\nu}\Phi(\bar{\textbf e}) \varepsilon^{-1}_r+\bar{\lambda}^{r-\nu}\Phi(\bar{\textbf e}) \varepsilon^{-1}_r+\bar{\lambda}^{r-\nu-1}\Phi(\bar{\textbf e}) \varepsilon^{-1}_r+\ldots +\bar{\lambda}\Phi(\bar{\textbf e}) \varepsilon^{-1}_r+\nonumber \\
&&~~~~~~~~~~~+\Phi(\bar{\textbf e}) \varepsilon^{-1}_r+a^{r-\nu}n)(a+1) \nonumber \\
&&~~~~~~~ \leq  \varepsilon^{-1}_r[\Phi(\bar{\textbf e})(r-\nu +2)+a^{-\nu}n](a+1). \nonumber
\end{eqnarray}	
Hence, we have
 \begin{displaymath}
\dim_B\mathcal{A}=\lim_{\varepsilon _r\rightarrow 0}\frac{\log N(\varepsilon _r)}{-\log \varepsilon _r}\leq 1+\frac{1}{r}\log_a{\Phi(\bar{\textbf e})(r-\nu +2) +a^{-\nu}n)(a+1)}.
\end{displaymath}
On the other hand $ \mathcal{A}$ is a curve in $\mathbf{R}^2$ and therefore $\dim_B\mathcal{A} \geq 1$. Hence $\dim_B\mathcal{A}=1$.     $\Box$\\ \\
\textbf{Remark 1}. In the case that $s_{i,k}(x)=s_{i,k}$(constant), if $\lambda=\bar{\lambda}=\underline{\lambda}>1$ , then $\dim_B\mathcal{ A}=1+\mathbf{\log}_a{\lambda}$. This is the estimation of box-counting dimension of RFISs in \cite{ban3}.

%
\section{Constructions of fractal surfaces using recurrent fractal interpolation curves}
In this section, we first introduce some constructions of fractal surfaces in which general fractal curves are used and compute their box-counting dimension. Next, this construction and recurrent fractal interpolation curves are combined.
\subsection{Construction of fractal surfaces by fractal curves and Lipschitz functions and their box-counting dimension}
We consider fractal surfaces only on the unit square $\left[0,1\right]\times\left[0,1\right]$ because our results can easily be generalized to any square $\left[a,b\right]\times\left[c,d\right]$. Let denote $I=\left[0,1\right]$ and $\mathrm{E}=I\times I$. \\ \indent
Let $f,~g:I\rightarrow\mathbf{R}$ be fractal curves (i.e. $f$ and $g$ are continuous and their graphs are fractal sets). Let $\lambda,~ \mu :I\rightarrow\mathbf{R}$ be continuous Lipschitz functions with Lipschitz constants $\mathrm{L}_\lambda$ and $~\mathrm{L}_\mu $ respectively. For $\left(x,y\right)\in\mathrm{E}$, we define a continuous function $F:\mathrm{E}\rightarrow\mathbf{R}$ by 
\begin{eqnarray}\label{F}
F\left(x,y\right)=\lambda\left(x\right)f\left(x\right)+\mu\left(y\right)g\left(y\right) .
\end{eqnarray}
The following theorem shows that the graph of this function $F$ is a fractal set.
\begin{Thot}\label{T:3}  
Let the function $F$ be given by (\ref{F}). Then the box-counting dimension $\mathrm{dim}_\mathrm{B}$ $Gr\left(F\right)$ of its graph is
\begin{eqnarray}\label{FF}
\mathrm{dim}_\mathrm{B}Gr\left(F\right)=1+\mathrm{Max}\left\{\mathrm{dim}_\mathrm{B}Gr\left(f\right),~\mathrm{dim}_\mathrm{B}Gr\left(g\right)\right\}\nonumber.
\end{eqnarray}
\end{Thot}

To prove the theorem, we need some basic lemmas on box-counting dimension. The following lemma is easily proved from the definition of the box-counting dimension.

\begin{lem}\label{L:2}   
Let $\mathcal{A}$ and $\mathcal{B}$ be fractal sets and $\mathrm{dim_{B}}\mathcal{A}\geq \mathrm{dim_{B}}\mathcal{B}$. Then
\begin{eqnarray}
&&\lim_{\delta \to 0}\frac{\log\left[\mathrm{N}_\delta\left(\mathcal{A}\right)+\mathrm{N}_\delta\left(\mathcal{B}\right)\right]}{-\log\delta}=\mathrm{dim}_\mathrm{B}\mathcal{A},\nonumber\\ 
&& \lim_{\delta \to 0}\frac{\log\left[\mathrm{N}_\delta\left(\mathcal{A}\right)+\mathrm{N}_\delta\left(\mathcal{B}\right)\pm \delta^{-1}\right]}{-\log\delta}=\mathrm{dim}_\mathrm{B}\mathcal{A},\nonumber\\
&&\lim_{\delta \to 0}\frac{\log\left[\mathrm{N}_\delta\left(\mathcal{A}\right)\pm \delta^{-1}\right]}{-\log\delta}=\mathrm{dim}_\mathrm{B}\mathcal{A}.\nonumber
\end{eqnarray}
\end{lem}

\begin{lem}\label{L:3}
Let the functions $f,~g $ be fractal curves. If we define a function $F':~\mathrm{E}\rightarrow\mathbf{R}$ by
\begin{displaymath}
F'\left(x,~y\right)=f\left(x\right)+g\left(y\right),\left(x,~y\right)\in\mathrm{E},
\end{displaymath} 
then the box-counting dimension $\mathrm{dim_{B}}Gr\left(F'\right)$ is as follows:
\begin{eqnarray}\label{F'}
\mathrm{dim}_\mathrm{B}Gr\left(F'\right)=1+\mathrm{Max}\left\{\mathrm{dim}_\mathrm{B}Gr\left(f\right),~\mathrm{dim}_\mathrm{B}Gr\left(g\right)\right\}\nonumber.
\end{eqnarray}
\end{lem}
\textbf{Proof.} Divide the interval $I$ into $n$ subintervals with the same length, denote the $i$th subinterval by $I_i=\left[\frac{i-1}{n},~\frac{i}{n}\right]$ and denote $\mathrm{E}_{ij}=I_i\times I_j$.
Let denote $N_\delta\left(\mathcal{A}\right)$ on a subinterval $I_i$ by $N^i_\delta(\mathcal{A})$ and on a subdomain $\mathrm{E}_{ij}$ by $N^{ij}_\delta(\mathcal{A})$. In the calculation of the box-counting dimension, we use $\frac{1}{n}$-mesh cubes for $n\in\mathbf{N} \left(n\geq 2\right)$ \\\indent
Evidently for any $i, ~j$ we have  $$\mathrm{R}_{F'}\left[\mathrm{E}_{ij}\right]=\mathrm{R}_f\left[\mathrm{I}_i\right]+\mathrm{R}_g\left[\mathrm{I}_j\right].$$ On the other hand, let $v:=R_{F'}[E_{ij}]\cdot (\frac{1}{n})^{-1}-[R_{F'}[E_{ij}]\cdot (\frac{1}{n})^{-1}]$. Here $[d]$ is the integer part of $d\in\mathbf{R}$ and $0\leq v < 1$. If $v=0$, then 
$$ R_{F'}[E_{ij}]\cdot \left(\frac{1}{n}\right)^{-1}=\left[R_{F'}[E_{ij}]\cdot \left(\frac{1}{n}\right)^{-1}\right]\leq N_{\frac{1}{n}}^{ij}(Gr(F'))\leq \left[R_{F'}[E_{ij}]\cdot \left(\frac{1}{n}\right)^{-1}\right]+1 $$
and if $0< v < 1$, then
$$\left[R_{F'}[E_{ij}]\cdot \left (\frac{1}{n}\right)^{-1}\right]+v\leq N_{\frac{1}{n}}^{ij}(Gr(F'))\leq \left[R_{F'}[E_{ij}]\cdot \left(\frac{1}{n}\right)^{-1}\right]+2$$
Similarly we have
$$\left[R_{f}[I_{i}]\cdot \left(\frac{1}{n}\right)^{-1}\right]\leq N_{\frac{1}{n}}^{i}(Gr(f))\leq \left[R_{f}[I_{i}]\cdot \left(\frac{1}{n}\right)^{-1}\right]+2,$$
$$\left[R_{g}[I_{j}]\cdot \left(\frac{1}{n}\right)^{-1}\right]\leq N_{\frac{1}{n}}^{j}(Gr(g))\leq \left[R_{g}[I_{j}]\cdot\left (\frac{1}{n}\right)^{-1}\right]+2.$$
On the other hand, noting that for any real numbers $m,n$ we have $[m]+[n]\leq[m+n]\leq[m]+[n]+1$, we have 
$$\left[R_{f}[I_{i}]\cdot \left(\frac{1}{n}\right)^{-1}\right]+\left[R_{g}[I_{j}]\cdot \left(\frac{1}{n}\right)^{-1}\right]\leq N_{\frac{1}{n}}^{ij}(Gr(F'))\leq \left[R_{f}[I_{i}]\cdot \left(\frac{1}{n}\right)^{-1}\right]+\left[R_{g}[I_{j}]\cdot \left(\frac{1}{n}\right)^{-1}\right]+3.$$
Therefore
\begin{eqnarray}
\mathrm{N}^i_{\frac{1}{n}}\left(Gr\left(f\right)\right)+\mathrm{N}^j_{\frac{1}{n}}\left(Gr\left(g\right)\right)-4&\leq&\mathrm{N}^{ij}_{\frac{1}{n}}\left(Gr\left(F'\right)\right)\leq \mathrm{N}^i_{\frac{1}{n}}\left(Gr\left(f\right)\right)+\mathrm{N}^j_{\frac{1}{n}}\left(Gr\left(g\right)\right)+3,  \nonumber\\
n\mathrm{N}_{\frac{1}{n}}\left(Gr\left(f\right)\right)+ n\mathrm{N}_{\frac{1}{n}}\left(Gr\left(g\right)\right)-4 n^2&\leq& \mathrm{N}_{\frac{1}{n}}\left(Gr\left(F'\right)\right)\leq n\mathrm{N}_{\frac{1}{n}}\left(Gr\left(f\right)\right)+ n\mathrm{N}_{\frac{1}{n}}\left(Gr\left(g\right)\right)+3 n^2. \nonumber 
\end{eqnarray}
From Lemma \ref{L:2} and the definition of the box-counting dimension we get the result. $\Box$\\\indent
\begin{lem}\label{L:4}
Let the functions $f,~\lambda $ be the same as the above ones. If we define the function $\lambda f:~I\rightarrow\mathbf{R}$  by $\left(\lambda f\right)\left(x\right)=\lambda\left(x\right)f\left(x\right), x\in I$, then the box-counting dimension $\mathrm{dim_{B}}Gr\left(\lambda f \right)$ is the same as that of the function $f$.
\end{lem}
\textbf{Proof.} For any $x,~x'\in I_i$,
\begin{eqnarray}
|\left(\lambda f\right)\left(x\right)-\left(\lambda f\right)\left(x'\right)|&=&|\lambda\left(x\right)f\left(x\right)-\lambda\left(x'\right)f\left(x'\right)|\nonumber\\
                                  &=&|\lambda\left(x\right)f\left(x\right)-\lambda\left(x\right)f\left(x'\right)+\lambda\left(x\right)f\left(x'\right)-
\lambda\left(x'\right)f\left(x'\right)| \nonumber\\
                                  &=&|\lambda\left(x\right)\left(f\left(x\right)-f\left(x'\right)\right)+\left(\lambda\left(x\right)-\lambda\left(x'\right)
\right)f\left(x'\right)|.\nonumber
\end{eqnarray}
Let denote $\mathrm{M}_f=\mathrm{Max}_{x\in\mathrm{I}}|f\left(x\right)|,~c_{\lambda_{1}}=\mathrm{Min}_{x\in\mathrm{I}}|
\lambda\left(x\right)|$ and $c_{\lambda_{2}}=\mathrm{Max}_{x\in\mathrm{I}}|\lambda\left(x\right)|$. We can suppose $c_{\lambda_1}>0$ (see Remark 2).
Then, for $n$ large enough,
\begin{eqnarray}
|\lambda\left(x\right)||f\left(x\right)-f\left(x'\right)|-\mathrm{L}_\lambda\frac{1}{n}|f\left(x'\right)|
 \leq  |\left(\lambda f\right)\left(x\right)-\left(\lambda f\right)\left(x'\right)|
\leq |\lambda\left(x\right)||f\left(x\right)-f\left(x'\right)|+\mathrm{L}_\lambda\frac{1}{n}|f\left(x'\right)|,\nonumber
\end{eqnarray}
\begin{eqnarray}
c_{\lambda_1}\mathrm{R}_f\left[{I}_i\right]-\mathrm{L}_\lambda\frac{1}{n}\mathrm{M}_f \leq  \mathrm{R}_{\lambda f}\left[{I}_i\right]\leq   c_{\lambda_2}\mathrm{R}_f\left[{I}_i\right]+\mathrm{L}_\lambda\frac{1}{n}\mathrm{M}_f.\nonumber
\end{eqnarray}
Thus
\begin{eqnarray}
c_{\lambda_1}{N}^i_{\frac{1}{n}}\left(Gr\left(f\right)\right)-m_1 &\leq &  {N}^i_{\frac{1}{n}}\left(Gr\left(\lambda f\right)\right)\leq   c_{\lambda_2}{N}^i_{\frac{1}{n}}\left(Gr\left(f\right)\right)+m_2,\nonumber \\
c_{\lambda_1}{N}_{\frac{1}{n}}\left(Gr\left(f\right)\right)-n\cdot m_1 &\leq &{N}_{\frac{1}{n}}\left(Gr\left(\lambda f\right)\right)\leq  c_{\lambda_2}\mathrm{N}_{\frac{1}{n}}\left(Gr\left(f\right)\right)+n\cdot m_2,\nonumber 
\end{eqnarray}
where $m_1=2c_{\lambda_1}+\mathrm{L}_\lambda\mathrm{M}_f$, $m_2=2+\mathrm{L}_\lambda\mathrm{M}_f$.
Taking logarithms, the definition of the box-counting dimension and Lemma \ref{L:2} give the result. $\Box$\\\\\indent
\textbf{Remark 2} In the case where $c_{\lambda_1}=0$, we choose $\lambda'\left(x\right)=\lambda\left(x\right)+c$ such that $c_{\lambda'_1}>0$ and define a function $\lambda' f: \mathrm{I}\rightarrow\mathbf{R}$ by $ \left(\lambda' f\right)\left(x\right)= \left(\lambda f\right)\left(x\right)+cf\left(x\right)$ for $x\in\mathrm{I}$. Then
$\left(\lambda f\right)\left(x\right)= \left(\lambda' f\right)\left(x\right)+\left(-cf\left(x\right)\right).$
Therefore, from Lemma \ref{L:2} and \ref{L:3}, we have $\mathrm{dim_B}Gr\left(\lambda f\right)=\mathrm{dim_B}Gr\left(f\right)$.\\\\
\textbf{Proof of Theorem 3.} Lemma \ref{L:3} and \ref{L:4} give the result of the theorem. $\Box$
\newtheorem{co}{Corollary}
\begin{co}\label{corollary1}
Let $f_i,~g_j$ be fractal curves and $\lambda_i,~\mu_j$ be (one-variable) Lipschitz functions for $i=1,\ldots,N,~j=1,\ldots,M$. Then the box-counting dimension $\mathrm{dim_B}Gr\left(F\right)$ of the function $F:\mathrm{E}\rightarrow\mathbf{R}$ defined for $\left(x,~y\right)\in\mathbf{E}$ by
\begin{displaymath}
F\left(x,~y\right)=\Sigma^N_{i=1}\lambda_i\left(x\right)f_i\left(x\right)+\Sigma^M_{j=1}\mu_j\left(y\right)g_j\left(y\right)\nonumber  
\end{displaymath} 
is as follows:
\begin{eqnarray}
\mathrm{dim_B}Gr\left(F\right)=1+\mathrm{Max}_{i,j}\left\{\mathrm{dim_B}Gr\left(f_i\right),~\mathrm{dim_B}Gr\left(g_j\right)\right\}.\nonumber
\end{eqnarray} 
\end{co}

The following theorem is proved analogously to Theorem \ref{T:3}. 
\begin{Thot}\label{gen}
Let the functions $f_i,~g_j$ be the same as in Theorem \ref{T:3} and functions $\lambda_i,~\mu_j~:\mathrm{E}\rightarrow\mathbf{R}$ continuous Lipschitz functions for $i=1,\ldots,N,~j=1,\ldots,M$. Then the function $F:\mathrm{E}\rightarrow\mathbf{R}$ defined for $\left(x,y\right)\in\mathrm{E}$ by
\begin{eqnarray}\label{ff}
F\left(x,y\right)=\Sigma^N_{i=1}\lambda_i\left(x,y\right)f_i\left(x\right)+\Sigma^M_{j=1}\mu_j\left(x,y\right)g_j\left(y\right)
\end{eqnarray} 
has the box-counting dimension of
\begin{eqnarray}  
\mathrm{dim_B}Gr\left(F\right)=1+\mathrm{Max}_{i,j}\left\{\mathrm{dim_B}Gr\left(f_i\right),~\mathrm{dim_B}Gr\left(g_j\right)\right\}.\nonumber
\end{eqnarray} 
\end{Thot}
\textbf{Remark 3} The fractal surfaces presented in the papers \cite{fal, WX, XFC, ZQL} are contained in the family of fractal surfaces defined by (\ref{F}) and (\ref{ff}).
%
\subsection{Construction of fractal surfaces by recurrent fractal interpolation curves}
As a simple application of the results of the preceding sections, we can easily construct fractal surfaces if we use at least 2 recurrent fractal interpolation curves. \\\indent
This methods can have wide applications in modeling  natural surfaces such as metals, rocks, terrains and so on or in interpolation of the data set in $\mathbf{R}$, $\mathbf{R}^2$. For one example, if data sets are given along the (x- and y-) boundaries (or several parallel lines to the boundaries) of a square, we can construct the fractal surfaces interpolating the data sets on the boundaries (or the parallel lines to the boundaries) by control of $\lambda_i(x,y)$ and $\mu _j(x,y)$ in  (\ref{ff}). \\\indent
The results of the calculations of the box-counting dimension in Section 4.1 show that the complexity of the fractal surfaces defined above are dominated by the fractal curves generating them. Thus, the more flexible a construction of fractal curves is, the more natural the fractal surface constructed by them is.  \\\indent
We can control the complexity and shape of the fractal surfaces constructed in Section 4.1 in our way by the vertical contractive function $s\left(x\right)$, the stochastic matrix $P$ and Lipschitz functions$g(x),~h(x)$ used in construction of recurrent fractal curves in section 3 and the Lipschitz functions $~\lambda$ and $\mu$ used in construction of fractal surfaces by combining fractal curves in section 4.1. So, the fractal surfaces presented in this paper could be more appropriate to model natural objects. \\\indent

\textbf{Example 2.} Construction of RFC and FS by RFC. \\\indent
The Fig. 1 shows the recurrent fractal curves generated from data set 
$ P=\{(1,20), $  $(0.25, 30), (0.5,$ $10), (0.75, 50), (1.0, 10)\} $ using the method of the Example 1 in Section 3 with different vertical scaling factor functions, where we use linear transformations $L_i$ and interpolation polynomials $g$ and $h$. The vertical scaling factors $\{s_1(x), s_2(x), s_3(x), s_4(x)\}$ are as follows: (a) $\{0.6,0.9,0.3,0.9\}$, (b) $\{-(x-1)^2/8+0.95, -(x-0.5)^2/7+0.95, 4(x-0.75)^2+0.6, 4(x-1)^2+0.7\}$, (c) $ \{0.9\sin(3.14x)$, $0.9\sin(3.14x),0.9\sin(3.14x),0.9\sin(3.14x)\}$, (d) $\{0.9\cos(x),0.98\cos(20x), 0.95 \cos(10x), 0.85\sin(50x)\}.$  \\\indent
The Fig. 2 shows the fractal surfaces constructed using recurrent fractal curves and different Lipschitz functions. The recurrent fractal curves $f(x)$, $g(y)$ are  (c), (d) in Fig. 1, respectively. In the Fig. 2, we used the Lipschitz functions $\lambda(x)=0.8$, $\mu(y)=0.5$ in (a) and $\lambda(x)= \cos(14\pi x)$, $\mu(y)=\sin(12\pi y)$ in (b), $\lambda(x,y)=(x - 0.5)^2+(y-0.5)^2$, $\mu(x,y)=-(x-0.5)^2-(y-0.5)^2$ in (c) and $\lambda(x,y)=(1-x)y$, $\mu(x,y)=x(1-y)$ in (d).  
\begin{figure}
\centering
\includegraphics[width=0.9\textwidth , height=9.5cm]{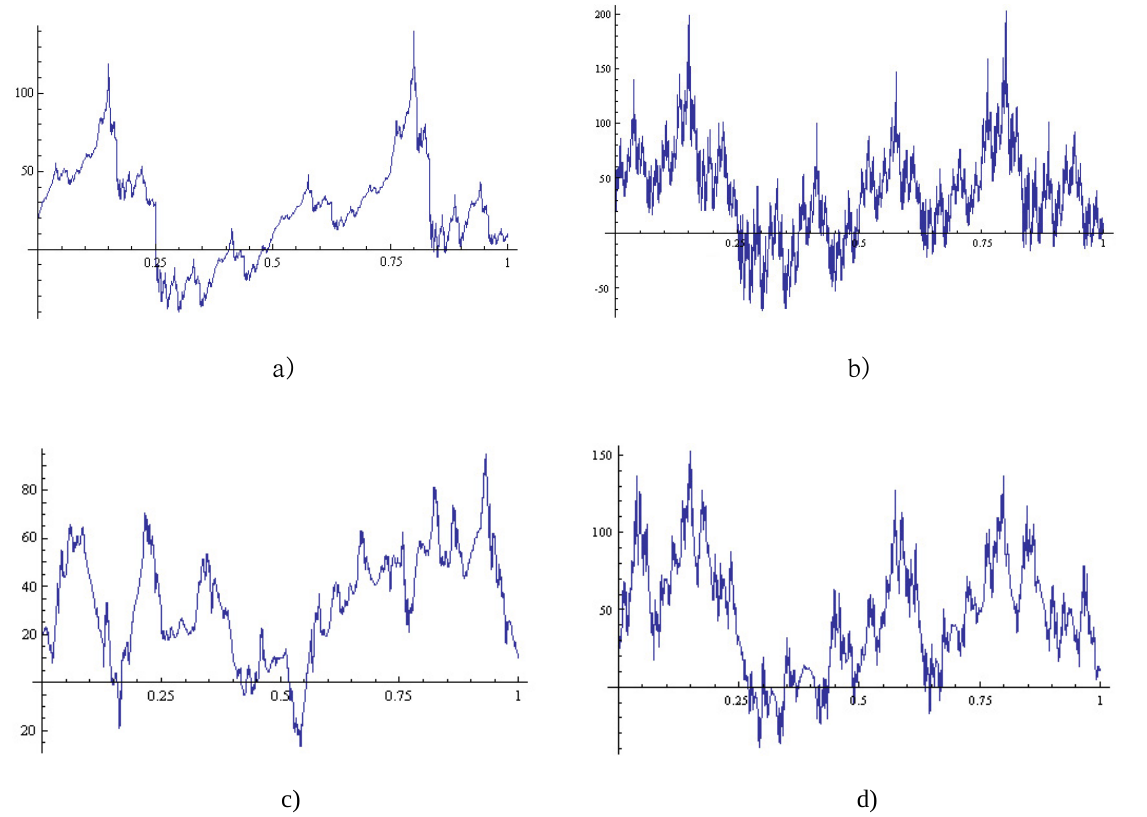}
\caption{ Recurrent fractal curves.}
\end{figure}
\begin{figure}
\centering
\includegraphics[width=0.9\textwidth , height=9.5cm]{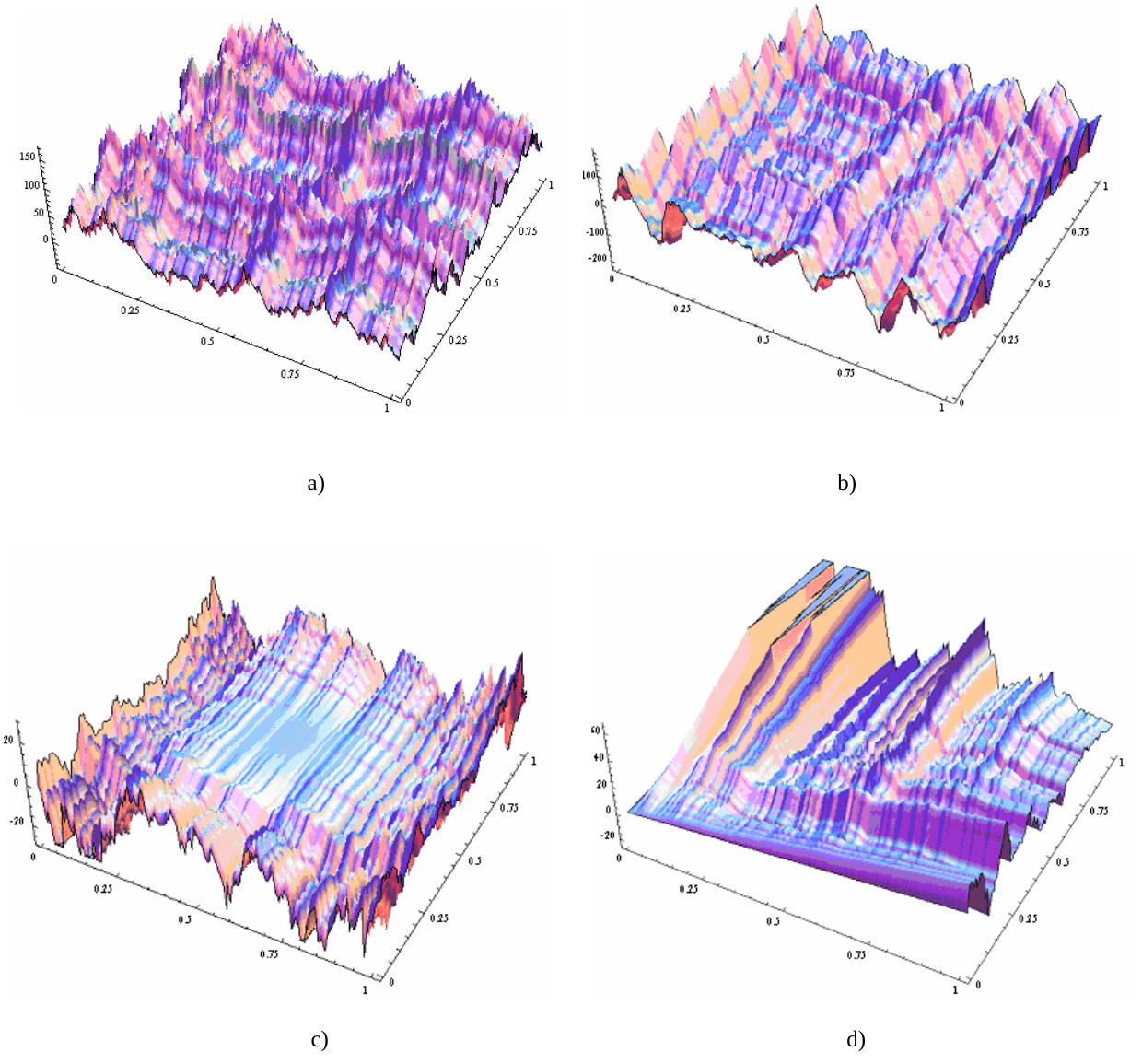}
\caption{ Fractal surfaces. }
\end{figure}


\end{document}